\documentclass[12pt,a4paper]{article}

\usepackage{amsmath,amsfonts, amssymb, graphicx}
\begin{document}
\bibliographystyle{abbrv}

\title{Log-concavity and the maximum entropy property of the Poisson distribution}
\author{Oliver Johnson\thanks{Statistical Laboratory, 
Centre for Mathematical Sciences, University of Cambridge, Wilberforce Rd, 
Cambridge, CB3 0WB, UK.
Email: {\tt otj1000@cam.ac.uk}\;
Fax: +44 1223 337956\; Phone: +44 1223 337946}}
\date{\today}
\maketitle

\newtheorem{theorem}{Theorem}[section]
\newtheorem{lemma}[theorem]{Lemma}
\newtheorem{proposition}[theorem]{Proposition}
\newtheorem{corollary}[theorem]{Corollary}
\newtheorem{conjecture}[theorem]{Conjecture}
\newtheorem{definition}[theorem]{Definition}
\newtheorem{example}[theorem]{Example}
\newtheorem{condition}{Condition}
\newtheorem{main}{Theorem}
\newtheorem{remark}[theorem]{Remark}
\hfuzz30pt

\def \outlineby #1#2#3{\vbox{\hrule\hbox{\vrule\kern #1%
\vbox{\kern #2 #3\kern #2}\kern #1\vrule}\hrule}}%
\def \endbox {\outlineby{4pt}{4pt}{}}%
\newenvironment{proof}
{\noindent{\bf Proof\ }}{{\hfill \endbox
}\par\vskip2\parsep}
\newenvironment{pfof}[2]{\removelastskip\vspace{6pt}\noindent
 {\it Proof  #1.}~\rm#2}{\par\vspace{6pt}}
\newenvironment{skproof}
{\noindent{\bf Sketch Proof\ }}{{\hfill \endbox
}\par\vskip2\parsep}

\newcommand{\var}{{\rm{Var\;}}}
\newcommand{\cov}{{\rm{Cov\;}}}
\newcommand{\tends}{\rightarrow \infty}
\newcommand{\C}[1]{{\bf ULC}(#1)}
\newcommand{\ep}{{\mathbb {E}}}
\newcommand{\pr}{{\mathbb {P}}}
\newcommand{\re}{{\mathbb {R}}}
\newcommand{\wt}[1]{\widetilde{#1}}
\newcommand{\wtt}[1]{\wt{\wt{#1}}}

\newcommand{\La}{\Lambda}
\newcommand{\I}{\mathbb {I}}
\newcommand{\Z}{{\mathbb {Z}}}
\newcommand{\bin}[2]{\binom{#1}{#2}}
\newcommand{\blah}[1]{}
\newcommand{\map}[1]{{\bf #1}}

\begin{abstract} \noindent
We prove that the Poisson distribution maximises entropy in the
class of ultra log-concave distributions, extending a result of Harremo\"{e}s.
The proof uses ideas concerning log-concavity, and a semigroup action involving
adding Poisson variables and thinning. We go on to show that the entropy
is a concave function along this semigroup.
\end{abstract}

\section{Maximum entropy distributions}
It is well-known that the distributions which maximise entropy under certain
very natural conditions take a simple form.
For example, among random variables with fixed mean and variance the
entropy is maximised by the normal distribution. Similarly, for random
variables with positive support and fixed mean, the entropy is maximised
by the exponential distribution.
The standard technique for proving such results uses the Gibbs inequality,
and establishes the fact that, given a function $R(\cdot)$ and fixing $\ep R(X)$, 
the maximum entropy density is of the
form $\alpha \exp(- \beta R(x))$ for constants $\alpha$ and $\beta$.
\begin{example} \label{ex:norm}
Fix mean $\mu$ and variance $\sigma^2$ and 
 write $\phi_{\mu,\sigma^2}$ for the 
density of $Z_{\mu,\sigma^2} \sim 
N(\mu, \sigma^2)$. For random variable $Y$ with
density $p_Y$ write
$\La(Y) = - \int p_Y(y) \log \phi_{\mu,\sigma^2}(y) dy$. Then
for any random variable $X$ with mean $\mu$, variance $\sigma^2$ and 
density $p_X$,
\begin{eqnarray}
\La(X) = - \int p_X(x) \log \phi_{\mu,\sigma^2}(x)  dx 
& = & 
\int p_X(x) \left( \frac{\log(2 \pi \sigma^2)}{2} 
+ \frac{(x-\mu)^2}{2\sigma^2} 
\right) dx \nonumber \\
& = & - \int \phi_{\mu,\sigma^2}(x) \log \phi_{\mu,\sigma^2}(x)  dx =
\La(Z_{\mu,\sigma^2}). 
\label{eq:mommatch}
\end{eqnarray}
This means that, for any random variable $X$ with mean $\mu$ and variance $\sigma^2$,
the entropy $H$ satisfies $H(X) \leq H(Z_{\mu,\sigma^2})$, since 
Equation (\ref{eq:mommatch}) gives that $\La(X) = \La(Z_{\mu,\sigma^2})
= H(Z_{\mu,\sigma^2})$,
\begin{eqnarray}
- H(X) + H(Z_{\mu,\sigma^2}) & = & \int p_X(x) \log p_X(x) dx - 
\int p_X(x) \log \phi_{\mu,\sigma^2}(x) dx. \label{eq:posi}
\end{eqnarray}
This expression is the relative entropy $D( X \| Z_{\mu,\sigma^2})$, and
is positive by the Gibbs inequality (see Equation (\ref{eq:logsum}) below),
with equality holding if and only if
 $p_X \equiv \phi_{\mu,\sigma^2}$.
\end{example} 
This maximum entropy result can be regarded as
the first stage in understanding the  Central Limit Theorem as a result concerning
maximum entropy. Note that both the class of variables 
with mean $\mu$ and variance $\sigma^2$ (over which the entropy is
maximised) and the maximum entropy variables $Z_{\mu,\sigma^2}$
are well-behaved on convolution. Further, the normalized
sum of IID copies of any random variable $X$ in this class converges in total
variation to the maximum entropy distribution $Z_{\mu,\sigma^2}$. The main
theorem of Barron \cite{barron} extends this to prove convergence
in relative entropy, assuming that $H(X) > - \infty$.

However, for functions $R$ 
where $\ep R(X)$ is not so well-behaved on convolution, 
the situation is more complicated. 
Examples of such random variables, for which we would hope to prove limit laws of a similar kind,  include the Poisson and Cauchy families. In particular, we 
would like to understand the ``Law of Small Numbers'' convergence to the 
Poisson distribution as a maximum entropy result.
Harremo\"{e}s proved in \cite{harremoes} that the Poisson random variables
$Z_\lambda$ (with mass function $\Pi_\lambda(x) = e^{-\lambda} \lambda^x/x!$
and mean $\lambda$) do satisfy a natural maximum entropy property.
\begin{definition} For each $\lambda \geq 0$ and $n \geq 1$ define the classes
$$ B_n(\lambda) = \biggl\{ S: \ep S = \lambda, S =
\sum_{i=1}^n X_i, \mbox{ where $X_i$ are independent Bernoulli variables} 
\biggr\}, $$
and $B_{\infty}(\lambda) = \bigcup_n B_n(\lambda)$. 
\end{definition}
\begin{theorem}[\cite{harremoes}, Theorem 8]
For each $\lambda \geq 0$, the entropy of any random variable in class
$B_{\infty}(\lambda)$ is less than or equal to the entropy of a 
Poisson random variable $Z_{\lambda}$:
$$ \sup_{S \in B_\infty(\lambda)} H(S) = H( Z_\lambda).$$
\end{theorem}
Note that Shepp and Olkin \cite{shepp} and Mateev \cite{mateev} also showed that the maximum entropy 
distribution in the class $B_n(\lambda)$ is Binomial($n,\lambda/n$).

In this paper, we show how this maximum entropy property 
relates to the property of log-concavity,
and give an alternative proof, which shows that $Z_\lambda$ is the maximum
entropy random variable in a larger class $\C{\lambda}$. 
\section{Log-concavity and main theorem}
First, recall the following definition:
\begin{definition} \label{def:logconc}
A non-negative 
sequence $(u(i), i \geq 0)$ is log-concave if, for all $i \geq 1$, 
\begin{equation} \label{eq:lcpropseq}
u(i)^2 \geq u(i+1) u(i-1).
\end{equation}
\end{definition}
We say that a random variable $V$ taking values in $\Z_+$ is log-concave if its
probability mass function $P_V(i) = \pr(V=i)$ forms a log-concave sequence.
Any random variable $S \in B_\infty$ is log-concave, which
is a corollary of the following theorem (see for example Theorem 1.2 on
P.394 of \cite{karlin3}).

\begin{theorem} \label{thm:hoggar} 
The convolution of any two log-concave sequences is log-concave.
\end{theorem}
Among random variables, the extreme cases of log-concavity are given
by the geometric family -- that is, geometric probability mass
functions are the only ones which achieve equality in Equation (\ref{eq:lcpropseq}) for all $i$. The argument of Example \ref{ex:norm} shows that discrete
entropy is maximised under a mean constraint by the geometric distribution.
Hence, in the class of log-concave random variables with a given mean, the geometric is both the extreme and the maximum entropy distribution.

Unfortunately, the sum of two geometric distributions
is a negative binomial distribution, which has a mass function which is
log-concave but no longer achieves equality in (\ref{eq:lcpropseq}). 
This means that under the condition of log-concavity
the extreme cases and the  maximum entropy family are not well-behaved under convolution. This 
suggests that log-concavity alone is too weak a condition to motivate an
entropy-theoretic understanding of the Law of Small Numbers.

A more restrictive condition than log-concavity is ultra log-concavity, defined
as follows:
\begin{definition} \label{def:ultralogconc}
A non-negative 
sequence $(u(i), i \geq 0)$ is ultra log-concave if the sequence $(u(i) i!, i \geq 0)$ is log-concave. That is, for all $i \geq 1$, 
\begin{equation} \label{eq:ulcpropseq}
i u(i)^2 \geq (i+1) u(i+1) u(i-1).
\end{equation}
\end{definition}

Note that in Pemantle \cite{pemantle}, Liggett \cite{liggett},
and Wang and Yeh \cite{wang3}, this property is
referred to as `ultra log-concavity of order $\infty$' -- see Equation 
(\ref{eq:ulcn}) below for the definition of ultra log-concavity of order
$n$.

An equivalent characterization of ultra log-concavity is that for any $\lambda$,
the sequence of ratios $(u(i)/\Pi_{\lambda}(i))$ is log-concave. This makes it
clear that among probability mass functions
the extreme cases of ultra log-concavity, in the sense of equality holding
in Equation (\ref{eq:ulcpropseq}) for each $i$, are exactly the Poisson 
family, which is preserved on convolution.
\begin{definition} \label{def:ulc} For any $\lambda \geq 0$,
define $\C{\lambda}$ to be the class of random variables $V$ with 
mean $\ep V = \lambda$ such that probability mass function $P_V$
is ultra log-concave, that is 
\begin{equation} \label{eq:ulc} i P_V(i)^2 \geq (i+1) P_V(i+1) P_V(i-1),
\mbox{ for all $i \geq 1$.}
\end{equation}
An equivalent characterization of the class $\C{\lambda}$ is that
the scaled score function introduced in \cite{johnson11} is
decreasing, that is
\begin{equation} \label{eq:ulc2} \rho_V(i) = \frac{(i+1) P_V(i+1)}
{\lambda P_V(i)} - 1
\mbox{ is a decreasing function in $i$.}
\end{equation}
\end{definition}
In Section \ref{sec:props} we discuss properties of the class $\C{\lambda}$.
For example, Lemma \ref{lem:ulcprops} shows
that (as for Harremo\"{e}s's $B_\infty(\lambda)$) the $\C{\lambda}$ are
well-behaved on convolution, and that $B_{\infty}(\lambda) \subset
\C{\lambda}$, with $Z_\lambda \in \C{\lambda}$. 

The main theorem of this paper is as follows:
\begin{theorem} \label{thm:main}
For any $\lambda \geq 0$, if $X \in \C{\lambda}$ then the entropy of $X$ satisfies
$$ H(X) \leq H(Z_{\lambda}),$$
with equality if and only if $X \sim Z_{\lambda}$.
\end{theorem}
We argue that this result gives the discrete analogue of the maximum entropy 
property of the normal distribution described in Example \ref{ex:norm}, 
since both the class $\C{\lambda}$ and
the family $Z_{\lambda}$ of maximum entropy random variables
are preserved on convolution, and since  $\C{\lambda}$ has another desirable property, that of ``accumulation''. That is, suppose
we fix $\lambda$ and take a triangular array of random variables $\{ X_i^{(n)} \}$, where for $i=1, \ldots, n$ the
$X_i^{(n)}$ are IID and in $\C{\lambda/n}$.
The techniques of \cite{johnson11} can be
extended to show that as $n \tends$ the sum 
$X_1^{(n)} + \ldots + X_n^{(n)}$ converges to $Z_{\lambda}$ in total
variation (and indeed in relative entropy).

It is natural to wonder whether Theorem \ref{thm:main} is optimal, or whether for 
each $\lambda$ there
exists a strictly larger class $C(\lambda)$ such that (i)  the $C(\lambda)$ are
well-behaved on convolution (ii) $Z_{\lambda}$ is the maximum entropy random
variable in each $C(\lambda)$ (iii) accumulation holds. We do not offer a
complete answer to this question
though, as discussed above, the class of log-concave variables
is too large and fails both conditions (i) and (ii).

For larger classes $C(\lambda)$, again consider a triangular array  where 
$\{ X_i^{(n)} \} \in C(\lambda/n)$. Write $p_n = \pr(X_i^{(n)} > 0)$
and $Q_n$ for the conditional
distribution $Q_n(x) = \pr(X_i^{(n)} = x | X_i^{(n)} > 0)$. If the classes
$C(\lambda)$ are large enough that we can find a subsequence $(n_k)$ such
that $Q_{n_k} \rightarrow Q$ and $\ep Q_{n_k} \rightarrow \ep Q$, then
the sum $X_1^{(n)} + \ldots + X_n^{(n)}$ converges to a compound Poisson
distribution $CP(\lambda/\ep Q, Q)$. Thus, if $C(\lambda)$ are large enough
that we can find a limit $Q \not\equiv \delta_1$ then the limit is not
equal to $Z_{\lambda}$ and so
the  property of accumulation fails. (Note that for $
X \in \C{\lambda}$ the $\pr(X \geq 2 | X > 0) \leq (\exp(\lambda) - \lambda - 1)/
\lambda$, so the only limiting conditional distribution is indeed $\delta_1$).

The proof of Theorem 
\ref{thm:main} is given in Sections \ref{sec:props}
 and \ref{sec:maxent}, and
is based on a family of maps $(\map{U}_\alpha)$ which we introduce in 
Definition \ref{def:umap} below. This map mimics
the role played by the Ornstein-Uhlenbeck semigroup in the normal case.
In the normal case, differentiating along this semigroup shows that the
probability densities satisfy a partial differential equation, the heat
equation, and hence that the derivative of relative entropy is the Fisher
information (a fact referred to as the de Bruijn identity -- see \cite{barron}).
This property is used by Stam \cite{stam} and Blachman
\cite{blachman} to prove the Entropy Power Inequality, which gives a sharp 
bound on the behaviour of continuous entropy on convolution. It is possible that 
a version of 
$\map{U}_{\alpha}$ may give a similar result for discrete entropy. 

As $\alpha$ varies between 1 and 0, the map
 $\map{U}_\alpha$ interpolates between a given random variable $X$ and a Poisson
random variable with the same mean. By establishing monotonicity properties 
with respect to $\alpha$, the maximum entropy result,
Theorem \ref{thm:main}, follows. The action of 
$\map{U}_{\alpha}$ is to thin $X$ and then to
add an independent Poisson random variable to it. 
In Section \ref{sec:maxent}, we use $\map{U}_{\alpha}$ 
to establish the maximum entropy
property of the Poisson distribution. The key expression is
Equation (\ref{eq:heateqn}), which shows that
the resulting probabilities satisfy an analogue of the heat equation.

We abuse terminology slightly in referring to $\map{U}_{\alpha}$ as a semigroup;
in fact (see Equation (\ref{eq:usemi}) below) $\map{U}_{\alpha_1} \circ \map{U}_{\alpha_2} = \map{U}_{\alpha_1 \alpha_2}$, so we would
require a reparametrization $\map{W}_{\theta} = \map{U}_{\exp(-\theta)}$
reminiscent of Bakry and \'{E}mery \cite{bakry}
to obtain the more familiar relation that
$\map{W}_{\theta_1} \circ \map{W}_{\theta_2} = 
\map{W}_{\theta_1 + \theta_2}$. However, in
Section \ref{sec:concave}, we argue that $\map{U}_{\alpha}$ has the 
`right' parametrization, by proving Theorem \ref{thm:main2} which shows that $H(\map{U}_{\alpha} X)$ is not only 
monotonically decreasing in $\alpha$, but is indeed a concave function of 
$\alpha$. We prove 
this by writing $H(\map{U}_{\alpha} X) = \La( \map{U}_{\alpha} X) - D( \map{U}_{\alpha} X \|
Z_{\lambda})$, and differentiating both terms. 

In contrast to conventions in Information Theory, throughout the paper
entropy is defined using logarithms to
base $e$. However, scaling by a factor of $\log 2$ restores
the standard definitions.
\section{Properties of $\C{\lambda}$ and definitions of maps} \label{sec:props}
In this section, we first note some results concerning properties of the classes $\C{\lambda}$, before defining actions of addition and thinning that will
be used to prove the main results of the paper.
\begin{lemma} \label{lem:ulcprops} For any $\lambda \geq 0$ and $\mu \geq 0$:
\begin{enumerate}
\item If $V \in \C{\lambda}$ then it is log-concave.
\item The Poisson random variable $Z_{\lambda} \in \C{\lambda}$.
\item The classes are closed on convolution: that is
for independent $U \in \C{\lambda}$ and $V \in \C{\mu}$, the sum $U+V \in
\C{\lambda + \mu}$.
\item $B_{\infty}(\lambda) \subset \C{\lambda}$.
\end{enumerate}
\end{lemma}
\begin{proof} Parts 1. and 2. follow from the definitions.
Theorem 1 of Walkup \cite{walkup} implies that Part 3. holds,
though a more direct proof is given by
Theorem 2 of Liggett \cite{liggett}.
Part 4. follows from Part 3., since any Bernoulli($p$) mass function scaled
by $\Pi_p$ is supported only on 2 points, so belongs to
$\C{p}$.
\end{proof}
We can give an alternative proof of Part 3 of Lemma \ref{lem:ulcprops}, using ideas of negative
association developed by Efron \cite{efron} and by Joag-Dev and Proschan
\cite{joag-dev}. The key result is
that if $U$ and $V$ are log-concave random variables, then for any decreasing
function $\phi$
$$ \ep [ \phi(U,V) | U+V = w] \mbox{ is a decreasing function of $w$.}$$
Now, the Lemma on P.471 of \cite{johnson11} shows that, writing 
$\alpha = \ep U/(\ep U + \ep V)$ and using the score function of Equation
(\ref{eq:ulc2}), for independent $U$ and $V$:
$$ \rho_{U+V}(w) = \ep[ \alpha \rho_U(U) + (1-\alpha) \rho_V(V) | U+V =w],$$
so that if $\rho_U$ and $\rho_V$ are decreasing, then so is $\rho_{U+V}$.

\begin{remark}
For each $n$, the Poisson mass function $\Pi_{\lambda}$ is
not supported on $[0,n]$ and hence $Z_\lambda 
\notin B_{n}(\lambda)$, so that $Z_\lambda \notin B_{\infty}(\lambda)$. Indeed,
we can see that the class of ultra log-concave random variables
is non-trivially
 larger than the class of Bernoulli sums. For all random variables $V \in B_n(\lambda)$, the Newton inequalities
(see for example Theorem 1.1 of Niculescu \cite{niculescu}) imply that
the scaled mass function $P_V(i)/\bin{n}{i}$ is log-concave, so that for all $i
\geq 1$:
\begin{equation} \label{eq:ulcn}
\frac{ i P_V(i)^2}{ (i+1) P_V(i+1) P_V(i+1)} \geq 
\frac{n-i+1}{n-i}.
\end{equation}
This is the property referred to by Pemantle \cite{pemantle} and Liggett
\cite{liggett} as ``ultra log-concavity of order $n$'', and is strictly more
restrictive than simply ultra log-concavity which (see Equation (\ref{eq:ulc})) 
only requires a lower bound of 1 on the right-hand side.
\end{remark}

Next we introduce the maps $\map{S}_{\beta}$ and $\map{T}_{\alpha}$  that will be key to our results.
\begin{definition} \label{def:stmap} Define the maps $\map{S}_\beta$ and
$\map{T}_{\alpha}$ which act as follows:
\begin{enumerate}
\item{For any $\beta \geq 0$, define the map $\map{S}_\beta$ that maps  
random variable $X$ to random variable 
$$\map{S}_{\beta} X \sim X + Z_{\beta},$$ 
where $Z_{\beta}$ is a Poisson$(\beta)$ random variable independent of $X$.}
\item{For any $0 \leq \alpha \leq 1$, define the map $\map{T}_\alpha$ that maps 
random variable $X$ to random variable 
$$\map{T}_{\alpha} X \sim \sum_{i=1}^X B_i(\alpha),$$
where $B_i(\alpha)$ are Bernoulli $(\alpha)$ random variables, independent of 
each other and of $X$.
This is the thinning operation introduced by R\'{e}nyi \cite{renyi4}.
}
\end{enumerate}
\end{definition}
We now show how these maps interact:
\begin{lemma} \label{lem:comm}
For any $0 \leq \alpha, \alpha_1, \alpha_2 \leq 1$ and for any $\beta, \beta_1, \beta_2 \geq 0$, the maps defined in Definition \ref{def:stmap} satisfy:
\begin{enumerate}
\item{$\map{S}_{\beta_1} \circ
\map{S}_{\beta_2} = \map{S}_{\beta_2} \circ \map{S}_{\beta_1} = \map{S}_{\beta_1 + \beta_2}.$}
\item{$\map{T}_{\alpha_1} \circ
\map{T}_{\alpha_2} = \map{T}_{\alpha_2} \circ \map{T}_{\alpha_1} = \map{T}_{\alpha_1 \alpha_2}.$}
\item{$\map{T}_{\alpha} \circ \map{S}_{\beta} = \map{S}_{\alpha \beta} \circ \map{T}_{\alpha}$.}
\end{enumerate}
\end{lemma}
\begin{proof} Part 1. follows immediately from the definition. To prove Part 2,
we write $B_i(\alpha_1 \alpha_2) = B_i(\alpha_1) B_i(\alpha_2)$ where 
$B_i(\alpha_1)$ and $B_i(\alpha_2)$ are independent, then for any $X$
$$ \map{T}_{\alpha_1 \alpha_2} X \sim
\sum_{i=1}^X B_i(\alpha_1) B_i(\alpha_2)
= \sum_{i: B_i(\alpha_1) = 1, i \leq X} B_i(\alpha_2)
= \sum_{i= 1}^{\map{T}_{\alpha_1} X} B_i(\alpha_2).$$ Part 3 uses the fact that
the sum of a Poisson number of IID Bernoulli random variables is itself Poisson.
This means that for any $X$
$$ \left( \map{T}_{\alpha} \circ \map{S}_{\beta}  \right) X \sim \sum_{i=1}^{\map{S}_\beta X} B_i(\alpha)
= \left( \sum_{i=1}^X B_i(\alpha) \right) + \left( \sum_{i=X+1}^{X+Z_\beta} 
B_i(\alpha) \right) \sim \map{T}_{\alpha} X + Z_{\alpha \beta} 
\sim \left( \map{S}_{\alpha \beta} \circ \map{T}_{\alpha} \right)X,$$
as required.
\end{proof}

\begin{definition} Define the two-parameter family of maps 
$$\map{V}_{\alpha,\beta} = \map{S}_{\beta} \circ \map{T}_{\alpha},
\mbox{ for $0 \leq \alpha \leq 1$, $\beta > 0$.} $$
\end{definition}
As in Stam \cite{stam} and Blachman \cite{blachman}, we will differentiate
along this family of maps, and see that the resulting probabilities satisfy a 
partial differential-difference equation.
\begin{proposition} \label{prop:heateqn}
Given $X$ with mean $\lambda$, 
writing $P_{\alpha}(z) = \pr(\map{V}_{\alpha,f(\alpha)} X = z)$, then
\begin{equation} \label{eq:heateqn}
 \frac{\partial }{\partial \alpha} P_{\alpha}(z)
= g(\alpha) (P_{\alpha}(z) - P_{\alpha}(z-1)) - \frac{1}{\alpha}
((z+1) P_{\alpha}(z+1) - 
z  P_{\alpha}(z)),\end{equation}
where $g(\alpha) = f(\alpha)/\alpha - f'(\alpha)$. Equivalently, 
$f(\alpha) = \alpha f(1) + \alpha \int_{\alpha}^1 g(\beta)/\beta d \beta$.
\end{proposition}
\begin{proof} 
We consider probability generating functions (pgfs). Notice that 
$$ \pr( \map{T}_{\alpha} X = z) = \sum_{x \geq z} \pr(X = x) 
\binom{x}{z} \alpha^z (1-\alpha)^{x-z},$$ so that
if $X$ has pgf $G_X(t) = \sum \pr(X = x) t^x$, then $\map{T}_{\alpha} X$ has pgf
$\sum_z t^z \sum_{x \geq z} \pr(X = x) 
\binom{x}{z} \alpha^z (1-\alpha)^{x-z} = 
\sum_x \pr(X = x) \sum_{z=0}^x \binom{x}{z} (t \alpha)^z (1-\alpha)^{x-z}
= G_X(t \alpha + 1-\alpha)$.

If $Y$ has pgf $G_Y(t)$ then $\map{S}_\beta Y$ has pgf
$G_Y(t) \exp( \beta(t-1))$. Overall then, $\map{V}_{\alpha,f(\alpha)} X$ has pgf 
\begin{equation} G_{\alpha}(t) =  
G_X( t \alpha + (1-\alpha)) \exp( f(\alpha) (t-1)),
\label{eq:pgf} \end{equation}
which satisfies
\begin{eqnarray*}
 \frac{\partial}{\partial \alpha} G_{\alpha}(t) = (1-t)
\left( \frac{1}{\alpha} \frac{\partial}{\partial t} G_{\alpha}(t)
- G_{\alpha}(t) g(\alpha) \right),
\end{eqnarray*}
and comparing coefficients the result follows.
\end{proof}
We now prove that both maps $\map{S}_{\beta}$ and $\map{T}_{\alpha}$ preserve ultra log-concavity.
\begin{proposition} \label{prop:upc}
If $X$ is an ultra log-concave random variable then for any $\alpha \in [0,1]$ and $\beta \geq 0$ random variables $\map{S}_{\beta} X$ and
$\map{T}_{\alpha} X$ are both ultra log-concave, and hence so is $\map{V}_{\alpha,\beta} X$.
\end{proposition}
\begin{proof} The first result follows by Part 3. of Lemma \ref{lem:ulcprops}.
We prove the second result using the case $f(\alpha) \equiv 0$
of Proposition \ref{prop:heateqn}, which tells us that writing
$P_{\alpha}(x) = \pr(\map{T}_{\alpha} X = x)$, the derivative
\begin{equation} \label{eq:pdesim}
 \frac{\partial }{\partial \alpha} P_{\alpha}(x)
= \frac{1}{\alpha} \left(x P_\alpha(x) - (x+1) P_{\alpha}(x+1) \right).
\end{equation}
Writing $g_\alpha(z) =
z P_{\alpha}(z)^2 - (z+1) P_{\alpha}(z+1) P_{\alpha}(z-1)$, Equation 
(\ref{eq:pdesim}) gives that for each $z$,
 \begin{eqnarray}
\frac{\partial }{\partial \alpha} g_{\alpha}(z) & = &
2z \frac{g_{\alpha}(z)}{\alpha} 
 + \frac{z+1}{\alpha} \big( (z+2) P_{\alpha}(z+2) P_{\alpha}(z-1) 
- z P_\alpha(z) P_{\alpha}(z+1) \big) \nonumber \\
& = &
\left( 2z - \frac{(z+2) P_{\alpha}(z+2)}{P_{\alpha}(z+1)} \right)
 \frac{g_{\alpha}(z)}{\alpha} - \frac{ z P_{\alpha}(z)}{\alpha P_{\alpha}(z+1)}
g_{\alpha}(z+1). \label{eq:change}
\end{eqnarray} 
We know that $P_\alpha$ is ultra log-concave for $\alpha =1$, and will show that this holds for smaller values of $\alpha$. 
Suppose that for some $\alpha$, $P_{\alpha}$ is ultra log-concave, so for each $z$,
$g_\alpha(z) \geq 0$. If for some $z$, $g_\alpha(z)
= 0$ then since $g_{\alpha}(z+1) \geq 0$,
Equation (\ref{eq:change}) simplifies to give 
$ \frac{\partial }{\partial \alpha} g_{\alpha}(z) \leq 0$. This means (by continuity) that there is no value of $z$ for which $g_{\alpha}(z)$ can become
negative as $\alpha$ gets smaller, so ultra log-concavity is preserved.
\end{proof}
\section{Maximum entropy result for the Poisson} \label{sec:maxent}
We now prove the maximum entropy property of the Poisson distribution within
the class $\C{\lambda}$. 
We choose a one-parameter family of maps $(\map{U}_{\alpha})$, which
have the property that they preserve the mean $\lambda$.
\begin{definition} \label{def:umap} Given mean $\lambda \geq 0$ and
$0 \leq \alpha \leq 1$, define
the combined map
$$ \map{U}_{\alpha} = \map{V}_{\alpha,\lambda(1-\alpha)}.$$
Equivalently $\map{U}_{\alpha} = \map{S}_{\lambda(1-\alpha)} \circ \map{T}_{\alpha} $
or $\map{U}_{\alpha}= \map{T}_{\alpha} \circ \map{S}_{\lambda(1/\alpha - 1)}.$
\end{definition}
Note that the maps $\map{U}_{\alpha}$ have  a 
semigroup-like structure -- by Lemma \ref{lem:comm} we know that
$(\map{S}_{\lambda (1-\alpha_1)} \circ \map{T}_{\alpha_1}) \circ
(\map{S}_{\lambda (1-\alpha_2)} \circ \map{T}_{\alpha_2}) 
= (\map{S}_{\lambda (1-\alpha_1)} \circ \map{S}_{\lambda \alpha_1 (1-\alpha_2)}) \circ
(\map{T}_{\alpha_1} \circ \map{T}_{\alpha_2}) = \map{S}_{\lambda (1 - \alpha_1 \alpha_2)}
\circ \map{T}_{\alpha_1 \alpha_2}$.
That is, we know that
\begin{equation} \label{eq:usemi}
\map{U}_{\alpha_1} \circ \map{U}_{\alpha_2} = \map{U}_{\alpha_1 \alpha_2}. \end{equation}

Equation (\ref{eq:heateqn}) can be simplified with the introduction of some
helpful notation. Define $\Delta$ and its
adjoint $\Delta^*$ by $\Delta p(x) = p(x+1) - p(x)$
and $\Delta^* q(x) = q(x-1) - q(x)$. These maps 
$\Delta$ and $\Delta^*$ are indeed
adjoint since for any functions $p,q$:
\begin{equation} \label{eq:adjo}
\sum_x \left( \Delta p(x) \right) q(x) = \sum_x (p(x+1) - p(x)) q(x) = \sum_x
p(x) (q(x-1) - q(x)) = \sum_x p(x) \left( \Delta^* q(x) \right). \end{equation}

We write $\rho_{\alpha}(z)$ for $\rho_{\map{U}_{\alpha} X}(z) = (z+1) P_{\alpha}(z+1)/\lambda P_{\alpha}(z) -1$.
Then, noting that $(z+1) P_{\alpha}(z+1)/\lambda - P_{\alpha}(z) = P_{\alpha}(z) 
\rho_{\alpha}(z) = \Pi_{\lambda}(z) \left( P_{\alpha}(z+1)/\Pi_{\lambda}(z+1)
- P_{\alpha}(z)/
\Pi_{\lambda}(z) \right)$, we can give two alternative reformulations of
Equation (\ref{eq:heateqn}) in the case where $\map{V}_{\alpha,f(\alpha)} = \map{U}_{\alpha}$.
\begin{corollary}
Writing $P_{\alpha}(z) = \pr(\map{U}_\alpha X = z)$:
\begin{equation} \label{eq:newheat}
\frac{\partial }{\partial \alpha} P_{\alpha}(z)  = 
\frac{\lambda}{\alpha} \Delta^* (P_{\alpha}(z) \rho_{\alpha}(z)).\end{equation}
Secondly, in a form more reminiscent of the heat equation:
$$
\frac{\partial }{\partial \alpha} P_{\alpha}(z)
= \frac{\lambda}{\alpha} \Delta^* \left( \Pi_{\lambda}(z) \Delta 
\left( \frac{P_{\alpha}(z)}{\Pi_{\lambda}(z)} \right) \right).$$
\end{corollary}
Note that we can also view $\map{U}_\alpha$ as the action of the M/M/$\infty$
queue. In particular Equation (\ref{eq:heateqn}), representing the 
evolution of probabilities under $\map{U}_{\alpha}$, is the adjoint of
$$Lf (z) = -\lambda \Delta \Delta^* f(z) + (z - \lambda) \Delta^* f(z),$$
representing the evolution of functions. This equation is the polarised
form of the infinitesimal generator of the M/M/$\infty$ queue, as described
in Section 1.1 of Chafa\"{i} \cite{chafai}. Chafa\"{i} uses this equation
to prove a number of inequalities concerning generalized entropy
functionals.

\begin{proof}{\bf of Theorem \ref{thm:main}}
Given random variable $X$ with mass function $P_X$, we 
define $\La(X) = - \sum_x P_X(x) \log \Pi_\lambda(x)$.
Notice that (as remarked by Tops{\o}e \cite{topsoe2}), the conditions required
in Example \ref{ex:norm} can be weakened. If $\La(X) \leq \La(Z_{\lambda}) =
H(Z_{\lambda})$ then adapting Equation (\ref{eq:posi}) gives that $-H(X) + 
H(Z_{\lambda}) \geq -H(X) + \La(X) =  D(X \| Z_{\lambda}) \geq 0$, and we can
deduce the maximum entropy property.

We will in fact show that if $X \in \C{\lambda}$ 
then $\La(\map{U}_\alpha X)$ is an decreasing function of
$\alpha$. In particular, since $\map{U}_0 X \sim Z_{\lambda}$, and $\map{U}_1 X 
\sim X$,
we deduce that $\La(X) \leq \La(Z_{\lambda})$. (A similar technique of
controlling the sign of the derivative is used by Blachman \cite{blachman}) 
and Stam \cite{stam} to prove the Entropy Power Inequality).

We simply differentiate and use Equations (\ref{eq:adjo}) and
(\ref{eq:newheat}).  Note that
\begin{eqnarray}
\frac{\partial}{\partial \alpha} \La(\map{U}_{\alpha} X)
& = & - \frac{\lambda}{\alpha} \sum_z
\Delta^* \big( P_{\alpha}(z) \rho_{\alpha}(z) \big) \log \Pi_{\lambda}(z) 
\nonumber \\
& = & - \frac{\lambda}{\alpha} \sum_z P_{\alpha}(z) \rho_{\alpha}(z)
\Delta \big( \log \Pi_{\lambda}(z) \big) \nonumber \\
& = & \frac{\lambda}{\alpha} \sum_z P_{\alpha}(z)
\rho_{\alpha}(z) \log \left( \frac{z+1}{\lambda}\right). \label{eq:id}
\end{eqnarray}
By assumption $X \in \C{\lambda}$, so by Proposition \ref{prop:upc}
$\map{U}_{\alpha} X \in \C{\lambda}$, which is equivalent to saying that the score function $\rho_{\alpha}(z)$ is decreasing in $z$. Further, note that 
$\sum_z P_{\alpha}(z) \rho_{\alpha}(z) = 0$.
Since $\log ((z+1)/\lambda)$ is increasing in $z$ (a fact which is equivalent
to saying that the Poisson mass function
$\Pi_{\lambda}(z)$ is itself log-concave),
$\frac{\partial}{\partial \alpha}
\La(\map{U}_{\alpha} X)$ is negative by 
Chebyshev's rearrangement lemma, since it is the covariance of
a decreasing and increasing function.

In fact, $\La(\map{U}_{\alpha} X)$ is strictly decreasing in 
$\alpha$, unless
$X$ is Poisson. This follows since equality holds in Equation (\ref{eq:id})
if and only if $\rho_{\alpha}(z) \equiv 0$, which characterizes the 
Poisson distribution.
\end{proof} 
\section{Concavity of entropy along the semigroup} \label{sec:concave}
In fact, rather than just showing that the Poisson distribution has a maximum 
entropy property, in this section we
establish a stronger result, as follows.
\begin{theorem}  \label{thm:main2}
If $X \in \C{\lambda}$, then the entropy of $\map{U}_{\alpha} X$ is a decreasing and concave function of $\alpha$, that is
$$ \frac{\partial}{\partial \alpha} H(\map{U}_{\alpha} X) \leq 0 
\hspace*{1cm} \mbox{ and } \hspace*{1cm}
\frac{\partial^2}{\partial \alpha^2} H(\map{U}_{\alpha} X) \leq 0,$$
with equality if and only if $X \sim \Pi_{\lambda}$.
\end{theorem}
\begin{proof} The proof is contained in the remainder of this section, and
involves
writing $H(\map{U}_{\alpha} X) = \La( \map{U}_{\alpha} X) - D( \map{U}_{\alpha} X \|
Z_{\lambda})$, and differentiating both terms. 

We have already shown in Equation (\ref{eq:id}) 
 that $\La(\map{U}_{\alpha} X)$ is decreasing in $\alpha$.
We show in Lemma \ref{lem:lconc} that it is concave in $\alpha$, 
and in Lemmas \ref{lem:dinc} and \ref{lem:dconv} respectively we show
that $D(\map{U}_{\alpha} X \| Z_{\lambda})$ is increasing and convex.
Some of the proofs of these lemmas are merely sketched, since they involve
long algebraic manipulations using Equation (\ref{eq:newheat}). \end{proof}

In the case of continuous random variables, Costa \cite{costa} uses the 
concavity of the entropy power on addition of an independent normal variable (a stronger result than concavity of entropy 
itself) to prove a version of the Entropy Power Inequality. We regard Theorem
\ref{thm:main2} as the first stage in a similar proof of a discrete form of the Entropy Power Inequality.
\begin{lemma} \label{lem:dinc}
For $X$ with mean $\lambda$, $D( \map{U}_{\alpha} X \|
Z_{\lambda})$ is an increasing function of $\alpha$.
\end{lemma}
\begin{proof}
We use Equation (\ref{eq:newheat}). Note that (omitting arguments for
the sake of brevity):
$$ \frac{\partial}{\partial \alpha} \sum
P_\alpha \log \left( \frac{P_{\alpha}}{\Pi_\lambda} \right)
=  \sum \frac{\partial P_{\alpha}}{\partial \alpha} 
\log \left( \frac{P_{\alpha}}{\Pi_\lambda} \right)
+ \sum \frac{\partial P_{\alpha}}{\partial \alpha} 
=  \sum \frac{\partial P_{\alpha}}{\partial \alpha} 
\log \left( \frac{P_{\alpha}}{\Pi_\lambda} \right).$$
This means that
\begin{eqnarray}
\frac{\partial}{\partial \alpha} D( \map{U}_{\alpha} X \| Z_{\lambda})
& = & \frac{\lambda}{\alpha} \sum_z
\Delta^* (P_{\alpha}(z) \rho_{\alpha}(z) )
\log \left( \frac{P_{\alpha}(z)}{\Pi_\lambda(z)} \right) \nonumber \\
& = & \frac{\lambda}{\alpha} \sum_z
P_{\alpha}(z) \rho_{\alpha}(z) 
\log \left( \frac{P_{\alpha}(z+1) \Pi_\lambda(z)} 
{P_{\alpha}(z) \Pi_\lambda(z+1)} \right) \nonumber \\
& = & \frac{\lambda}{\alpha} \sum_z
P_{\alpha}(z) \rho_{\alpha}(z) 
\log \left( 1 + \rho_{\alpha}(z) \right). \label{eq:ddiff}
\end{eqnarray}
Now, as in \cite{johnson11}, we write $\wt{P}_{\alpha}(z) = (z+1) 
P_{\alpha}(z+1)/\lambda$. $\wt{P}_{\alpha}$ is often referred to as the
size-biased version of $P_{\alpha}$, and is a probability mass function
because $\map{U}_\alpha$ fixes the mean. Notice that $\rho_{\alpha}(z)
= \wt{P}_{\alpha}(z)/P_{\alpha}(z) - 1$, so that we can rewrite
Equation (\ref{eq:ddiff}) as 
\begin{eqnarray}
\frac{\lambda}{\alpha} \sum_z (\wt{P}_{\alpha}(z) - P_{\alpha}(z))
\log \left( \frac{ \wt{P}_{\alpha}(z)}{P_{\alpha}(z)} \right) 
\label{eq:ddiff2} 
& = &  \frac{\lambda}{\alpha} \left(
D( P_{\alpha} \| \wt{P}_{\alpha} ) + D( \wt{P}_{\alpha} \| P_{\alpha}) 
\right) \geq 0. 
\end{eqnarray}
This quantity is a
symmetrised version of the relative entropy, and was originally
introduced by Kullback and Leibler in \cite{kullback2}. 
\end{proof}
\begin{lemma} 
Using the definitions above, \label{lem:lconc} if $X \in \C{\lambda}$ then
$\La(\map{U}_{\alpha} X)$ is a concave
function of $\alpha$. It is strictly concave unless $X$ is Poisson.
\end{lemma}
\begin{skproof}
Using Equations (\ref{eq:newheat}) and (\ref{eq:id}), it can be shown that
\begin{eqnarray*}
\frac{\partial^2}{\partial \alpha^2} \La(\map{U}_{\alpha} X)
\blah{ & = & \frac{\partial}{\partial \alpha} \left(
\frac{\lambda}{\alpha} \sum_z \left( \frac{(z+1) P_{\alpha}(z+1)}
{\lambda} - P_{\alpha}(z) \right)
\log \left( \frac{z+1}{\lambda}\right) \right) \\
& = & -\frac{\lambda}{\alpha^2} \sum_z \left( 
\frac{(z+1) P_{\alpha}(z+1)}{\lambda} - P_{\alpha}(z) \right)
\log \left( \frac{z+1}{\lambda}\right)  \\
& & + \sum_z \frac{\lambda^2}{\alpha^2} 
\left( \frac{(z+1)}{\lambda}
\Delta^* (P_{\alpha}(z+1) \rho_{\alpha}(z+1) )
 - \Delta^* (P_{\alpha}(z) \rho_{\alpha}(z) ) \right)
\log \left( \frac{z+1}{\lambda}\right)  \\
& = &
\frac{\lambda}{\alpha^2} \sum_z 
\Delta^* \biggl( (z+1) P_{\alpha}(z+1) \rho_{\alpha}(z+1) \biggr)
\log \left( \frac{z+1}{\lambda}\right)  \\
& & - \frac{\lambda^2}{\alpha^2} \sum_{z} 
\Delta^* \biggl(P_{\alpha}(z) \rho_{\alpha}(z) \biggr) 
\log \left( \frac{z+1}{\lambda}\right)  \\
& = &
\frac{\lambda}{\alpha^2} \sum_z 
(z+1) P_{\alpha}(z+1) \rho_{\alpha}(z+1) 
\log \left( \frac{z+2}{z+1} \right)  \\
& & - \frac{\lambda^2}{\alpha^2} \sum_{z} P_{\alpha}(z) \rho_{\alpha}(z)  
\log \left( \frac{z+2}{z+1}\right)  \\ }
& = &
\frac{\lambda^2}{\alpha^2} \sum_z 
P_{\alpha}(z) \rho_{\alpha}(z) 
\left( \frac{z}{\lambda} \log \left( \frac{z+1}{z} \right)
- \log \left( \frac{z+2}{z+1}\right)  \right).
\end{eqnarray*}
Now, the result follows in the same way as before, since for any $\lambda$ the
function 
$z/\lambda \log( (z+1)/z) - \log( (z+2)/(z+1))$ is increasing,
so $\frac{\partial^2}{\partial \alpha^2} \La(\map{U}_{\alpha} X) \geq 0$.
\end{skproof}
Taking a further derivative of Equation (\ref{eq:ddiff2}), we can show that
(the proof is omitted for the sake of brevity):
\begin{lemma} \label{lem:dpde}
The relative entropy $D(\map{U}_{\alpha} X \| Z_{\lambda})$ satisfies
\begin{eqnarray*}
\lefteqn{\frac{\partial^2}{\partial \alpha^2} \sum_z
P_\alpha(z) \log \left( \frac{P_{\alpha}(z)}{\Pi_\lambda(z)} \right) } \\
& = & \frac{\lambda^2}{\alpha^2}  
\sum_z (\wtt{P}_{\alpha}(z) - 2\wt{P}_{\alpha}(z)
+ P_{\alpha}(z) ) \log \left( \frac{ \wtt{P}_{\alpha}(z) P_{\alpha}(z)}
{\wt{P}_{\alpha}(z)^2} \right) + 
\sum_z P_{\alpha}(z) \left( \frac{1}{P_{\alpha}(z)} 
\frac{ \partial P_{\alpha}}{\partial \alpha}(z) \right)^2
\end{eqnarray*}
where $\wt{P}_{\alpha}(z) = (z+1) 
P_{\alpha}(z+1)/\lambda$ and $\wtt{P}_{\alpha}(z)
= (z+2)(z+1) P_{\alpha}(z+2)/\lambda^2$.
\end{lemma} \blah{
\begin{proof}
Again, we use Equation (\ref{eq:newheat}). Note that
$\frac{\partial}{\partial \alpha} \wt{P}_{\alpha}(z)
= (z+1) \Delta^*( P_{\alpha}(z+1) \rho_{\alpha}(z+1))/\alpha$. Taking a
further derivative of Equation (\ref{eq:ddiff2}), we deduce that
\begin{eqnarray*}
\lefteqn{\frac{\partial^2}{\partial \alpha^2} \sum_z
P_\alpha(z) \log \left( \frac{P_{\alpha}(z)}{\Pi_\lambda(z)} \right) } \\
& = & 
\frac{\lambda}{\alpha^2} \sum_z \Delta^* \biggl( (z+1) P_{\alpha}(z+1) 
\rho_{\alpha}(z+1) - \lambda P_{\alpha}(z) \rho_{\alpha}(z) \biggr)
\log \left( \frac{ \wt{P}_{\alpha}(z)}{P_{\alpha}(z)} \right) \\
& & + \frac{\lambda}{\alpha} \sum_z (\wt{P}_{\alpha}(z) - P_{\alpha}(z))
\left( \frac{\partial \wt{P}_{\alpha}}{\partial \alpha}(z)
\frac{1}{\wt{P}_{\alpha}(z)}
- \frac{\partial P_{\alpha}}{\partial \alpha}(z)
\frac{1}{P_{\alpha}(z)} \right)  \\
& = & 
\frac{\lambda^2}{\alpha^2} \sum_z \left( \frac{(z+2)(z+1) P_{\alpha}(z+2)}{
\lambda^2} - 2 \frac{(z+1) P_{\alpha}(z+1)}{\lambda} + 1 \right)
\log \left( \frac{ (z+2) P_{\alpha}(z+2) P_{\alpha}(z)}{
(z+1) P_{\alpha}^2(z+1)} \right) \\
& & + \frac{\lambda}{\alpha} \sum_z P_{\alpha}(z) \rho_\alpha(z)
\left( \frac{\partial P_{\alpha}}{\partial \alpha}(z+1)
\frac{1}{P_{\alpha}(z+1)}
- \frac{\partial P_{\alpha}}{\partial \alpha}(z)
\frac{1}{P_{\alpha}(z)} \right)  \\
& = & \frac{\lambda^2}{\alpha^2} 
\sum_z (\wtt{P}_{\alpha}(z) - 2\wt{P}_{\alpha}(z)
+ P_{\alpha}(z) ) \log \left( \frac{ \wtt{P}_{\alpha}(z) P_{\alpha}(z)}
{\wt{P}_{\alpha}(z)^2} \right)  \\
& & + \frac{\lambda}{\alpha} \sum_z \Delta^* (P_{\alpha}(z) \rho_\alpha(z))
\frac{\partial P_{\alpha}}{\partial \alpha}(z)
\frac{1}{P_{\alpha}(z)}  \\
& = & \frac{\lambda^2}{\alpha^2}  
\sum_z (\wtt{P}_{\alpha}(z) - 2\wt{P}_{\alpha}(z)
+ P_{\alpha}(z) ) \log \left( \frac{ \wtt{P}_{\alpha}(z) P_{\alpha}(z)}
{\wt{P}_{\alpha}(z)^2} \right) + 
\sum_z P_{\alpha}(z) \left( \frac{\partial P_{\alpha}}{\partial \alpha}(z)
\frac{1}{P_{\alpha}(z)} \right)^2,
\end{eqnarray*}
as required, using Equation (\ref{eq:newheat}).
\end{proof} }
\begin{lemma} \label{lem:dconv}
For $X$ with mean $\lambda$ and $\var X \leq \lambda$, 
$D( \map{U}_{\alpha} X \|
Z_{\lambda})$ is a convex function of $\alpha$. It is a strictly convex
function unless $X$ is Poisson.
\end{lemma}
\begin{proof}
Notice that the map $\map{T}_{\alpha}$ scales
the $r$th falling moment of $X$ by $\alpha^r$.
This means that $\var \map{T}_{\alpha} X = \alpha^2 \var X + \alpha(1-\alpha) \lambda$, so that $\var \map{U}_{\alpha} X = \alpha^2 \var X + \lambda(1-\alpha^2)$. Hence, the condition $\var X \leq \lambda$ implies that
for all $\alpha$, $\var \map{U}_{\alpha} X \leq \lambda$. Equivalently,
$S := \sum_z \wtt{P}_{\alpha}(z)
= \ep (\map{U}_\alpha X)(\map{U}_\alpha X -1)/\lambda^2 < 1$, 

We will use the log-sum inequality, which
is equivalent to the Gibbs inequality, and states that for positive sequences
$(a_i)$ and $(b_i)$ (not necessarily summing to 1),
\begin{equation} \label{eq:logsum}
D(a_i \| b_i) =
\sum_i a_i \log(a_i/b_i) \geq \left( \sum_i a_i \right) \log \left( 
\frac{\sum_i a_i}{\sum_i b_i} \right).
\end{equation}
Since $\log u \leq u -1$, this simplifies further to give
$D(a_i \| b_i) \geq (\sum_i a_i) \left( \log (\sum_i a_i) + 1 - \sum_i b_i
\right)$.

We express the first term of Lemma 
\ref{lem:dpde} as a sum of relative entropies, and recall that 
$\sum_z P_{\alpha}(z) = 1$ and $\sum_z \wt{P}_{\alpha}(z) = 1$, simplifying
the bounds on the second and third terms:
\begin{eqnarray}
\lefteqn{
\frac{\lambda^2}{\alpha^2} \left(
D \left( \wtt{P}_{\alpha} \left\|
\frac{\wt{P}_{\alpha}^2}{P_{\alpha}} \right) \right.
+ 2 D \left( \wt{P}_{\alpha} \left\|
\frac{\wtt{P}_{\alpha} P_{\alpha}}{\wt{P}_{\alpha}} \right) \right.
+ D \left( P_{\alpha} \left\|
 \frac{\wt{P}^2}{\wtt{P}} \right) \right. \right)} \nonumber \\
\blah{ & = & \frac{\lambda^2}{\alpha^2} \left(
S \log \left( \frac{S}
{\sum_z \wt{P}_{\alpha}(z)^2/P_{\alpha}(z)} \right)
- 2 \log \left( \sum_z \frac{\wtt{P}_{\alpha}(z) P_{\alpha}(z)}{
\wt{P}_{\alpha}(z)} \right)
- \log \left( \sum_z \frac{\wt{P}_{\alpha}(z)^2}{\wtt{P}_{\alpha}(z)} 
\right) \right) \nonumber \\
& \geq & \frac{\lambda^2}{\alpha^2} \left( S \log S -
S \left( \sum_z \frac{\wt{P}_{\alpha}(z)^2}{P_{\alpha}(z)} - 1 \right)
- 2 \left( \sum_z \frac{\wtt{P}_{\alpha}(z) P_{\alpha}(z)}{
\wt{P}_{\alpha}(z)} - 1 \right)
- \left( \sum_z \frac{\wt{P}_{\alpha}(z)^2}{\wtt{P}_{\alpha}(z)} 
-1 \right) \right) \nonumber \\ }
& \geq & \frac{\lambda^2}{\alpha^2} \left( S \log S + S -
S \sum_z \frac{(z+1)^2 P_{\alpha}(z+1)^2}{\lambda^2 P_{\alpha}(z)} 
+ 2 - 2 
\sum_z \frac{(z+1) P_{\alpha}(z+1) P_{\alpha}(z-1)}{\lambda P_{\alpha}(z)}
 \right. \nonumber \\
& & \left. \hspace*{1cm} + 1 - 
\sum_z \frac{(z-1) P_{\alpha}(z-1)^2}{z P_{\alpha}(z)}  \right).
\label{eq:terma}
\end{eqnarray}
Using Equation (\ref{eq:heateqn}) we can expand the second (Fisher) term 
of Lemma \ref{lem:dpde} as 
\begin{eqnarray}
\blah{ \lefteqn{ \frac{1}{\alpha^2}
\left( \lambda^2
\sum_z \frac{P_{\alpha}(z-1)^2}{P_{\alpha}(z)}
+ \sum_z \frac{ (z+1)^2 P_{\alpha}(z+1)^2}{ P_{\alpha}(z)} 
+ 2 \lambda \sum_z (z+1) \frac{P_{\alpha}(z-1) P_{\alpha}(z+1)}{
P_{\alpha}(z)} \right. } \nonumber \\ 
& & \left. + \sum_z P_\alpha(z) (z+ \lambda)^2
- 2 \lambda \sum_z P_{\alpha}(z) (\lambda + z+1)
- 2 \lambda \sum_z (z+1) P_{\alpha}(z+1)( \lambda +z)  \right)  \nonumber \\ }
& = & \frac{\lambda^2}{\alpha^2} \left( - 3  - \frac{\ep 
(\map U_{\alpha} X)^2}{\lambda^2} +
 \sum_z \frac{ (z+1)^2 P_{\alpha}(z+1)^2}{
\lambda^2 P_{\alpha}(z)} + 2  \sum_z \frac{(z+1) P_{\alpha}(z+1) 
P_{\alpha}(z-1)}{\lambda P_{\alpha}(z)} \right. \nonumber \\
& & \hspace*{1cm} \left. + \sum_z \frac{P_{\alpha}(z-1)^2}{P_{\alpha}(z)}
 \right). \label{eq:termb}
\end{eqnarray}
Adding Equations (\ref{eq:terma}) and (\ref{eq:termb}), and since
$S = \ep (\map U_{\alpha} X)^2/\lambda^2 - 1/\lambda$, 
we deduce that
\begin{eqnarray}
 \frac{ \partial^2}{\partial \alpha^2} D(\map{U}_{\alpha} X \| \Pi_{\lambda})
 \geq \frac{\lambda^2}{\alpha^2}
\left( S \log S + 
(1-S) \sum_z \frac{ (z+1)^2 P_{\alpha}(z+1)^2}{\lambda^2 P_{\alpha}(z)} 
+ \sum_z \frac{P_{\alpha}(z-1)^2}{z P_{\alpha}(z)} 
-\frac{1}{\lambda}  \right). \label{eq:combined}
\end{eqnarray}
Finally we exploit Cram\'{e}r-Rao type relations which bound the two
remaining quadratic terms from below. Firstly, as in \cite{johnson11}:
\begin{equation} \label{eq:cr1} 0 \leq
\sum_z P_{\alpha}(z) \left( \frac{ (z+1) P_{\alpha}(z+1)}{ \lambda P_{\alpha}(z)}
-1 \right)^2
= \sum_z \frac{ (z+1)^2 P_{\alpha}(z+1)^2}{\lambda^2 P_{\alpha}(z)} - 1.
\end{equation}
Similarly, a weighted version of the Fisher information term of Johnstone
and MacGibbon \cite{johnstone} gives that:
\begin{equation} \label{eq:cr2} 0 \leq
\sum_z P_{\alpha}(z) z \left( \frac{ P_{\alpha}(z-1)}{ z P_{\alpha}(z)}
- \frac{1}{\lambda} \right)^2
= \sum_z \frac{ P_{\alpha}(z-1)^2}{z P_{\alpha}(z)} - 
\frac{1}{\lambda}.\end{equation}
(Note that in Equations (\ref{eq:cr1}) and (\ref{eq:cr2}), equality holds if and only if $P_{\alpha} \equiv \Pi_{\lambda}$).
Substituting Equations (\ref{eq:cr1}) and (\ref{eq:cr2})
 in Equation (\ref{eq:combined}), we deduce that
\begin{eqnarray*}
 \frac{ \partial^2}{\partial \alpha^2} D(\map{U}_{\alpha} X \| Z_{\lambda})
 \geq \frac{\lambda^2}{\alpha^2}
\left( S \log S + 1-S \right) \geq 0,
\end{eqnarray*}
with equality if and only if $S = 1$.
\end{proof}
Combining these lemmas, the proof of Theorem \ref{thm:main2} is complete, since
ultra log-concavity of $X$ implies that $\var X \leq \ep X$, as
$ \sum_x P_X(x) x ( (x+1) P_X(x+1)/P_X(x) - \lambda) \leq 0$ since it is again the covariance of an increasing and decreasing function.
\section*{Acknowledgment}
The author would like to thank Christophe Vignat, Ioannis Kontoyiannis,
Peter Harremo\"{e}s, Andrew Barron and Mokshay Madiman for useful discussions 
concerning this paper, and would like to thank
EPFL Lausanne and Yale University for financial support on visits to 
these colleagues. The author would also like to thank Djalil Chafa\"{i}
for explaining the connection with the M/M/$\infty$ queue, and two anonymous
referees for their very helpful comments, including a simplified proof of 
Proposition \ref{prop:upc}.


\end{document}